\documentstyle[amssymb,12pt]{article}

\begin{document}

\title{{\bf Poisson-gradient dynamical systems with convex potential}}
\author{Iulian Duca, Ana-Maria Teleman and Constantin Udri\c ste}
\date{}
\maketitle

\begin{abstract}
The basic aim is to extend some results and concepts of non-autonomous
second order differential systems with convex potentials to the new context
of multi-time Poisson-gradient PDE systems with convex potential. In this
sense, we prove that minimizers of a suitable action functional are multiple
periodical solutions of a Dirichlet problem associated to the Euler-Lagrange
equations. Automatically, these are solutions of the associated multi-time
Hamiltonian equations.
\end{abstract}



\section{Poisson-gradient PDEs}

There are two methods to study the periodic solutions of boundary problems
attached to some partial derivative equations (PDEs):

- the method of Fourier expansions in terms of eigenfunctions of a PDE
operator (the method of separation of variables);

- the method of minimizers of suitable action functionals.

Our paper refers to the second method, continuing the ideas in the papers
[11], [14], [19]. We start with the set $T_{0}=\left[ 0,T^{1}\right] \times
...\times \left[ 0,T^{p}\right] \subset R^{p}$ determined by the diagonal
points $O=\left( 0,...,0\right) ,T=\left( T^{1},...,T^{p}\right) ,$\ and
with the Sobolev space $W_{T}^{1,2}$ of the functions $u\in L^{2}\left[
T_{0},R^{n}\right] $, having weak derivatives $\displaystyle\frac{\partial u%
}{\partial t}\in L^{2}\left[ T_{0},R^{n}\right] .$ The weak derivatives are
defined using the space $C_{T}^{\infty }$ of all indefinitely differentiable
multiple T-periodic\ functions from $R^{p}$ into $R^{n}$.

We denote by $H_{T}^{1}$ the Hilbert space associated to the Sobolev space $%
W_{T}^{1,2}.$ The euclidean structure on $H_{T}^{1}$ is given by the scalar
product 
\[
\langle u,v\rangle =\int_{T_{0}}\left( \delta _{ij}u^{i}\left( t\right)
v^{j}\left( t\right) +\delta _{ij}\delta ^{\alpha \beta }\frac{\partial u^{i}%
}{\partial t^{\alpha }}\left( t\right) \frac{\partial v^{j}}{\partial
t^{\beta }}\left( t\right) \right) dt^{1}\wedge ...\wedge dt^{p} 
\]
and the associated Euclidean norm. These are induced by the scalar product
(Riemannian metric) 
\[
G=\left( 
\begin{array}{cc}
\delta _{ij} & 0 \\ 
0 & \delta ^{\alpha \beta }\delta _{ij}
\end{array}
\right) 
\]
on $R^{n+np}$ (see the jet space $J^{1}\left( T_{0},R^{n}\right) $).

Let $t=\left( t^{1},...,t^{p}\right) $ be a generic point in $R^{p}$. Then
the opposite faces of the parallelepiped $T_{0}$ can be described by the
equations 
\[
S_{i}^{-}:t^{i}=0,S_{i}^{+}:t^{i}=T^{i} 
\]
for each $i=1,...,p.$

Suppose the function $u\left( t\right) $ has a weak Laplacian $\Delta u$ and 
$u\rightarrow F\left( t,u\right) $ is a convex function. In these
hypothesis, we formulate some conditions in which the Dirichlet problem
(associated to a Poisson-gradient PDE system) 
$$
\Delta u\left( t\right) =\nabla F\left( t,u\left( t\right) \right) \eqno(1) 
$$
$$
u\mid _{S_{i}^{-}}=u\mid _{S_{i}^{+}},\frac{\partial u}{\partial t}\mid
_{S_{i}^{-}}=\frac{\partial u}{\partial t}\mid _{S_{i}^{+}},i=1,...,p\eqno%
(2) 
$$
has solution. To do that, we denote 
\[
\left| \frac{\partial u}{\partial t}\right| ^{2}=\delta ^{\alpha \beta
}\delta _{ij}\frac{\partial u^{i}}{\partial t^{\alpha }}\frac{\partial u^{j}%
}{\partial t^{\beta }} 
\]
and we use the Lagrangian 
$$
L\left( t,u\left( t\right) ,\frac{\partial u}{\partial t}\right) =\frac{1}{2}%
\left| \frac{\partial u}{\partial t}\right| ^{2}+F\left( t,u\left( t\right)
\right) \eqno
(3) 
$$
and the action 
$$
\varphi \left( u\right) =\int_{T_{0}}L\left( t,u\left( t\right) ,\frac{%
\partial u}{\partial t}\right) dt^{1}\wedge ...\wedge dt^{p}\eqno
(4) 
$$
Then, using minimizing sequences, we show that the action $\varphi $ has a
minimum point $u$ (extremal, solution of the Poisson-gradient dynamical
system (1), satisfying the boundary conditions (2)). Consequently the
solution $u$ is multiple periodical, with the reduced period $T=\left(
T^{1},...,T^{p}\right) $. Our arguments extend those of the book [6],
Theorems 1.4, 1.5, 1.7 and 1.8, which are dedicated to single-time problems.

\section{Periodic \ solutions \ of \ Poisson-gradient PDEs}

Let us show that some conditions upon the potential $F$ ensure periodic
solutions for the problem (1)+(2).

{\bf Theorem 1} {\it Let $F:T_{0}\times R^{n}\rightarrow R,(t,x)\rightarrow
F(t,x)$ and }$\left| x\right| {\it =}\sqrt{\delta _{ij}x^{i}x^{j}}${\it . We
consider that }$F\left( t,x\right) $ {\it is measurable in $t$ for any $x\in
R^{n}$ and of class $C^{1}$ in $x$ for any} $t\in T_{0}.$

If {\it there exist $a\in C^1\left( R^{+},R^{+}\right) $ with the derivative $a'$ bounded from above\
and $b\in C\left( T_{0},R^{+}\right) $ such that 
\[
\left| F\left( t,x\right) \right| \leq a\left( \left| x\right| \right)
b\left( t\right) ,\left| \nabla _{x}F\left( t,x\right) \right| \leq a\left(
\left| x\right| \right) b\left( t\right) , 
\]
for any $x\in R^{n}$ and any $t$ $\in T_{0}$, then the action (4) is of
class $C^{1}$. }

{\bf Proof}. The reasons are similar to those in [15, Theorem 3].

{\bf Corollary 2} {\it The some hypothesis as in Theorem 1. If }${\it u\in
H_{T}^{1}}${\it \ is a solution of the equation }${\it \varphi ^{\prime
}\left( u\right) =0}${\it \ (critical point), then the function }${\it u}$%
{\it \ has a weak Laplacian }${\it \bigtriangleup u}${\it \ (the Jacobian
matrix }${\it \displaystyle
\frac{\partial u}{\partial t}}$ {\it has a weak divergence) and } 
\[
{\it \bigtriangleup u=\nabla F(t,u(t))}
\]
{\it a.e. on }${\it T_{0}}${\it \ and } 
$$
{\it u\mid _{S_{i}^{-}}=u\mid _{S_{i}^{+}},\frac{\partial u}{\partial t}\mid
_{S_{i}^{-}}=\frac{\partial u}{\partial t}\mid _{S_{i}^{+}}.}\eqno
(5)
$$
{\bf Proof.} We build the function 
\[
\Phi :[-1,1]\rightarrow R,
\]
\[
\Phi \left( \lambda \right) =\varphi \left( u+\lambda v\right) =
\]
\[
\int_{T_{0}}\left[ \frac{1}{2}\left| \frac{\partial }{\partial t}\left(
u\left( t\right) +\lambda v\left( t\right) \right) \right| ^{2}+F\left(
t,u\left( t\right) +\lambda v\left( t\right) \right) \right] dt^{1}\wedge
...\wedge dt^{p},
\]
where $v\in C_{T}^{\infty }.$ The point $\lambda =0$ is a critical point of $%
\Phi $ if and only if the point $u$ is a critical point of $\varphi $.
Consequently 
\[
0=\left\langle \varphi ^{\prime }\left( u\right) ,v\right\rangle
=\int_{T_{0}}\left[ \delta ^{\alpha \beta }\delta _{ij}\frac{\partial u^{i}}{%
\partial t^{\alpha }}\frac{\partial v^{j}}{\partial t^{\beta }}+\delta
_{ij}\nabla ^{i}F\left( t,u\left( t\right) \right) v^{j}\left( t\right) %
\right] dt^{1}\wedge ...\wedge dt^{p},
\]
for all $v\in H_{T}^{1}$ and hence for all $v\in C_{T}^{\infty }.$ Using the definition of the weak divergence, 
\[
\int_{T_{0}}\delta ^{\alpha \beta }\delta _{ij}\frac{\partial u^{i}}{%
\partial t^{\alpha }}\frac{\partial v^{j}}{\partial t^{\beta }}dt^{1}\wedge
...\wedge dt^{p}
=-\int_{T_{0}}\delta ^{\alpha \beta }\delta _{ij}\frac{\partial ^{2}u^{i}}{%
\partial t^{\alpha }\partial t^{\beta }}v^{j}dt^{1}\wedge ...\wedge dt^{p},
\]
the Jacobian matrix $\displaystyle
\frac{\partial u}{\partial t}$ has weak divergence (the function $u$ has a
weak Laplacian) and 
\[
\bigtriangleup u\left( t\right) =\nabla F\left( t,u\left( t\right) \right) 
\]

a.e. on $T_{0}$. Also, the existence of the weak derivatives $\displaystyle
\frac{\partial u}{\partial t}$ and weak divergence $\bigtriangleup u$
implies that 
\[
u\mid _{S_{i}^{-}}=u\mid _{S_{i}^{+}},\frac{\partial u}{\partial t}\mid
_{S_{i}^{-}}=\frac{\partial u}{\partial t}\mid _{S_{i}^{+}}. 
\]

{\bf Remark}. If the function $u$ is at least of class $C^2$, then the definition of the weak divergence of the Jacobian matrix ${\frac{\partial u}{\partial t}}$
(or of the weak Laplacian $\bigtriangleup u$) coincides with the classical definition. This fact is obvious if we have in mind the formula of 
{\it integration by parts}

\[
\int_{T_{0}}\delta ^{\alpha \beta }\delta _{ij}\frac{\partial u^{i}}{%
\partial t^{\alpha }}\frac{\partial v^{j}}{\partial t^{\beta }}dt^{1}\wedge
...\wedge dt^{p} 
\]
\[
=\int_{T_{0}}\delta ^{\alpha \beta }\delta _{ij}\frac{\partial }{\partial
t^{\alpha }}\left( \frac{\partial u^{i}}{\partial t^{\alpha }}v^{j}\right)
dt^{1}\wedge ...\wedge dt^{p}-\int_{T_{0}}\delta ^{\alpha \beta }\delta _{ij}%
\frac{\partial ^{2}u^{i}}{\partial t^{\alpha }\partial t^{\beta }}%
v^{j}dt^{1}\wedge ...\wedge dt^{p}. 
\]

{\bf Corollary 3} {\it The some hypothesis as in Theorem 1. If }${\it \left|
x\right| \rightarrow =}${\it \ implies }${\it \displaystyle
\int_{T_{0}}F\left( t,x\right) dt^{1}\wedge ...\wedge dt^{p}\rightarrow
\infty }${\it \ and }${\it F\left( t,x\right) }${\it \ is convex in }${\it x}
${\it \ for any }${\it t\in T_{0},}${\it \ then there exists a function }$%
{\it u}${\it \ that is a solution of the boundary value problem (5).}

{\bf Proof. \ }Let $G:R^{n}\rightarrow R,\,G\left( x\right) =\displaystyle
\int_{T_{0}}F\left( t,x\right) dt^{1}\wedge ...\wedge dt^{p}.$ By
assumptions, the convex function $G$ has a minimum point $x=\overline{x}$.
Consequently, $\nabla G(\overline{x})=\displaystyle\int_{T_{0}}\nabla
F\left( t,\overline{x}\right) dt^{1}\wedge ...\wedge dt^{p}=0$.

Let $\left( u_{k}\right) $ be a minimizing sequence for the action (4). We
use the decomposition $u_{k}=\overline{u}_{k}+\widetilde{u}_{k}$, where $%
\overline{u}_{k}=\displaystyle\int_{T_{0}}u_{k}\left( t\right) dt^{1}\wedge
...\wedge dt^{p}$. The convexity of $F$ implies 
\[
F\left( t,u_{k}\left( t\right) \right) \geq F\left( t,\overline{x}\right)
+\left( \nabla F\left( t,\overline{x}\right) ,u_{k}\left( t\right) -%
\overline{x}\right) . 
\]
It follows 
\[
\varphi \left( u_{k}\right) \geq \frac{1}{2}\int_{T_{0}}\left| \frac{%
\partial u_{k}}{\partial t}\right| ^{2}dt^{1}\wedge ...\wedge dt^{p}+ 
\]
\[
+\int_{T_{0}}F\left( t,\overline{x}\right) dt^{1}\wedge ...\wedge
dt^{p}+\int_{T_{0}}\left( \nabla F\left( t,\overline{x}\right) ,u_{k}\left(
t\right) -\overline{x}\right) dt^{1}\wedge ...\wedge dt^{p} 
\]
\[
=\frac{1}{2}\int_{T_{0}}\left| \frac{\partial u_{k}}{\partial t}\right|
^{2}dt^{1}\wedge ...\wedge dt^{p}+\int_{T_{0}}F\left( t,\overline{x}\right)
dt^{1}\wedge ...\wedge dt^{p}+ 
\]
\[
+\int_{T_{0}}\left( \nabla F\left( t,\overline{x}\right) ,\widetilde{u}%
_{k}\left( t\right) \right) dt^{1}\wedge ...\wedge dt^{p}. 
\]
On the other hand, by Schwartz inequality, we can write 
\[
\left( \nabla F\left( t,\overline{x}\right) ,\widetilde{u}_{k}\left(
t\right) \right) \leq \left| \nabla F\left( t,\overline{u}\right) \right|
\left| \widetilde{u}_{k}\left( t\right) \right| \leq a\left( \left| 
\overline{x}\right| \right) b\left( t\right) \left| \widetilde{u}_{k}\left(
t\right) \right| . 
\]
Consequently, 
\[
\varphi \left( u_{k}\right) \geq \frac{1}{2}\int_{T_{0}}\left| \frac{%
\partial u_{k}}{\partial t}\right| ^{2}dt^{1}\wedge ...\wedge dt^{p}+ 
\]
\[
+\int_{T_{0}}F\left( t,\overline{x}\right) dt^{1}\wedge ...\wedge
dt^{p}-a\left( \left| \overline{x}\right| \right) \int_{T_{0}}b\left(
t\right) \widetilde{u}_{k}\left( t\right) dt^{1}\wedge ...\wedge dt^{p} 
\]
\[
\geq \frac{1}{2}\int_{T_{0}}\left| \frac{\partial u_{k}}{\partial t}\right|
^{2}dt^{1}\wedge ...\wedge dt^{p}+ 
\]
\[
+\int_{T_{0}}F\left( t,\overline{x}\right) dt^{1}\wedge ...\wedge
dt^{p}-a\left( \left| \overline{x}\right| \right) b_{0}\int_{T_{0}}\left| 
\widetilde{u}_{k}\left( t\right) \right| dt^{1}\wedge ...\wedge dt^{p}, 
\]
where $b_{0}=\displaystyle\max_{t\in T_{0}}b\left( t\right) $. Using the
Wirtinger inequality for multiple integral, we find 
\[
\varphi \left( u_{k}\right) \geq \frac{1}{2}\int_{T_{0}}\left| \frac{%
\partial u_{k}}{\partial t}\right| ^{2}dt^{1}\wedge ...\wedge dt^{p}+ 
\]
\[
+\int_{T_{0}}F\left( t,\overline{x}\right) dt^{1}\wedge ...\wedge
dt^{p}-a\left( \left| \overline{x}\right| \right) b_{0}C_{1}\left(
\int_{T_{0}}\left| \frac{\partial u_{k}}{\partial t}\right| ^{2}dt^{1}\wedge
...\wedge dt^{p}\right) ^{\frac{1}{2}}, 
\]
with $C_{1}>0$. Thus 
\[
\varphi \left( u_{k}\right) \geq \frac{1}{2}\int_{T_{0}}\left| \frac{%
\partial u_{k}}{\partial t}\right| ^{2}dt^{1}\wedge ...\wedge
dt^{p}+C_{2}-C_{3}\left( \int_{T_{0}}\left| \frac{\partial u_{k}}{\partial t}%
\right| ^{2}dt^{1}\wedge ...\wedge dt^{p}\right) ^{\frac{1}{2}}, 
\]
and, consequently, the function of degree two in the right hand member must
be a decreasing restriction, i.e., there exists $C_{4}>0$, such that \newline
$\displaystyle\int_{T_{0}}\left| \displaystyle\frac{\partial u_{k}}{\partial
t}\right| ^{2}dt^{1}\wedge ...\wedge dt^{p}<C_{4}$. It follows 
\[
\int_{T_{0}}\left| \frac{\partial \widetilde{u}_{k}}{\partial t}\right|
^{2}dt^{1}\wedge ...\wedge dt^{p}<C_{4} 
\]
and so $\left\| \widetilde{u}_{k}\right\| <C_{5}$.

Again, the convexity of $F$ leads to 
\[
F\left( t,\frac{\overline{u}_{k}}{2}\right) =F\left( t,\frac{1}{2}\left(
u_{k}\left( t\right) -\widetilde{u}_{k}\left( t\right) \right) \right) \leq 
\frac{1}{2}F\left( t,u_{k}\left( t\right) \right) +\frac{1}{2}F\left( t,-%
\widetilde{u}_{k}\left( t\right) \right) , 
\]
$\forall t\in T_{0}$, $\forall k\in N$, so 
\[
F\left( t,u_{k}\left( t\right) \right) \geq 2F\left( t,\frac{\overline{u}%
_{k}\left( t\right) }{2}\right) -F\left( t,-\widetilde{u}_{k}\left( t\right)
\right) . 
\]
Consequently 
\[
\varphi \left( u_{k}\right) =\frac{1}{2}\int_{T_{0}}\left| \frac{\partial
u_{k}}{\partial t}\right| ^{2}dt^{1}\wedge ...\wedge
dt^{p}+\int_{T_{0}}F\left( t,u_{k}\left( t\right) \right) dt^{1}\wedge
...\wedge dt^{p}\geq 
\]
\[
\geq \frac{1}{2}\int_{T_{0}}\left| \frac{\partial u_{k}}{\partial t}\right|
^{2}dt^{1}\wedge ...\wedge dt^{p}+, 
\]
\[
+2\int_{T_{0}}F\left( t,\frac{\overline{u}_{k}}{2}\right) dt^{1}\wedge
...\wedge dt^{p}-\int_{T_{0}}F\left( t,-\widetilde{u}_{k}\left( t\right)
\right) dt^{1}\wedge ...\wedge dt^{p} 
\]
and hence $\varphi \left( u_{k}\right) \geq 2\displaystyle%
\int_{T_{0}}F\left( t,\displaystyle\frac{\overline{u}_{k}}{2}\right)
dt^{1}\wedge ...\wedge dt^{p}-c_{6}$.

This means that $||\overline{u}_{k}||\!\!\nrightarrow \;\infty $. So the
sequence $\left( \overline{u}_{k}\right) $ is bounded and implicitly the
sequence $\left( u_{k}\right) $ is bounded in $H_{T}^{1}$. The Hilbert space 
$H_{T}^{1}$ is reflexive. By consequence, the sequence $\left( u_{k}\right) $
(or one of his subsequence) is weakly convergent in $H_{T}^{1}$ with the
limit $u$. The Mazur's theorem assure that there exists a sequence $\left(
v_{k}\right) $ with the general term $\ v_{k}=\displaystyle%
\sum_{j=1}^{k}\alpha _{kj}u_{j},\displaystyle\sum_{j=1}^{k}\alpha
_{kj}=1,\alpha _{kj}\geq 0$ , which is strongly converges to $u$ in $%
H_{ST}^{1}$.

Now we consider $c>$\underline{$\lim $} $\varphi \left( u_{k}\right) $.
Going if necessary to a subsequence, we can assume that $c>\varphi \left(
u_{k}\right) $ for all $k\in N^{\ast }.$ Since $\varphi $ is lower
semi-continuous in $H_{T}^{1}$ and $\varphi $ is convex, we obtain 
\[
\varphi \left( u\right) \leq \underline{\lim }\varphi \left( v_{k}\right)
\leq \underline{\lim }\left( \displaystyle\sum_{j=1}^{k}\alpha _{kj}\varphi
\left( u_{j}\right) \right) \leq \left( \displaystyle\sum_{j=1}^{k}\alpha
_{kj}\right) c=c. 
\]

Because $c>$\underline{$\lim $ }$\varphi \left( u_{k}\right) $ is \
arbitrary, we have $\varphi \left( u\right) \leq $\underline{$\lim $ }$%
\varphi \left( u_{k}\right) .$

Thus, the action $\varphi \left( u\right) $ has a minimum point $u$ in $%
H_{T}^{1},$ and so $u$ is a solution of the problem (5).

Thanks to the properties of the strictly convex functions, we can reinforce
the previous theorem. For that, we recall two equivalent properties of a
strictly convex function $G\in C^{1}\left( R^{n},R\right) $:

1) $G$ has a critical point $\overline{x}\in R^{n};$

2) $G\left( x\right) \rightarrow \infty $ when $\left| x\right| \rightarrow
\infty .$

{\bf Theorem 2}. {\it We consider the problem (1)+(2).\ Suppose $%
F:T_{0}\times R^{n}\rightarrow R,(t,x)\rightarrow F(t,x)$ has the properties}%
:

1) {\it $F\left( t,x\right) $} {\it is measurable in $t$ for any $x\in R^{n}$
and of class $C^{1}$ in $x$ for any} ${\it t\in T}_{{\it 0}}{\it .}$

2) {\it There exist $a\in C^1\left( R^{+},R^{+}\right) $ with the derivative $a'$ bounded from above
and $b\in C\left( T_{0},R^{+}\right) $ so that} 
\[
{\it \left| F\left( t,x\right) \right| \leq a\left( \left| x\right| \right)
b\left( t\right) ,\left| \nabla _{x}F\left( t,x\right) \right| \leq a\left(
\left| x\right| \right) b\left( t\right) ,} 
\]
{\it for any $x\in R^{n}$ and any $t$ $\in $}${\it T_{0}.}$

3) {\it The function $F\left( t,\cdot \right) $ is strictly convex for any $%
t\in T$}${\it _{0}.}$

{\it Then, the following statements are equivalent}:

1) {\it The problem $\left( 1\right) +\left( 2\right) $ has solutions;}

2) {\it There exists $\overline{x}\in R^{n}$ so that $\displaystyle%
\int_{T_{0}}\nabla F\left( t,\overline{x}\right) dt^{1}\wedge ...\wedge
dt^{p}=$}${\it 0;}$

3) {\it $\displaystyle\int_{T_{0}}F\left( t,x\right) dt^{1}\wedge ...\wedge
dt^{p}\rightarrow \infty $ as $\left| x\right| \rightarrow \infty $.}

{\bf Proof} (see single-time case in [6, Theorem 1.8]). Let us prove that 1)
implies 2):

We suppose that $u\left( t\right) $ is a solution of the problem (1)+(2). By
integration we obtain 
\[
\displaystyle\sum_{i=1}^{p}\displaystyle\int_{T_{0}}\displaystyle\frac{%
\partial ^{2}u^{j}}{\partial t^{i^{2}}}dt^{1}\wedge ...\wedge
dt^{p}=\int_{T_{0}}\frac{\partial F}{\partial u^{j}}\left( t,u\left(
t\right) \right) dt^{1}\wedge ...\wedge dt^{p}. 
\]
From the boundary conditions it results 
$$
\displaystyle\int_{T_{0}}\nabla F\left( t,u\left( t\right) \right)
dt^{1}\wedge ...\wedge dt^{p}=0.\eqno
(6) 
$$
On the other hand, the function $G\left( x\right) =\displaystyle
\int_{T_{0}}F(t,x)dt^{1}\wedge ...\wedge dt^{p}$ is strictly convex, because
the function $F(t,\cdot )$ is strictly convex.

We suppose $u=\widetilde{u}+\overline{u}$, $\overline{u}=\displaystyle%
\int_{T_{0}}u\left( t\right) dt^{1}\wedge ...\wedge dt^{p},$ 
\[
\widetilde{G}\left( x\right) =\displaystyle\int_{T_{0}}F\left( t,x+u\left(
t\right) \right) dt^{1}\wedge ...\wedge dt^{p}. 
\]
From (6) we have $\nabla \widetilde{G}\left( \overline{u}\right) =0$. From
the properties of a strictly convex function, mentioned above, $\widetilde{G}%
\left( x\right) $ tends to $\infty $ when $\left| x\right| $ tends to $%
\infty $. Because the function $F\left( t,\cdot \right) $ is strictly
convex, we obtain: 
\[
\widetilde{G}\left( x\right) \leq \frac{1}{2}\int_{T_{0}}F\left( t,2x\right)
dt^{1}\wedge ...\wedge dt^{p}+\frac{1}{2}\int_{T_{0}}F\left( t,2u\left(
t\right) \right) dt^{1}\wedge ...\wedge dt^{p}=\frac{1}{2}G\left( 2x\right)
+C. 
\]
For $\left| x\right| \rightarrow \infty $, $\widetilde{G}\left( x\right)
\rightarrow \infty $ \ and consequently $G\left( 2x\right) \rightarrow
\infty $ and $G\left( x\right) \rightarrow \infty $. According to the
properties of $G$, there exists $\overline{x}$ so that $\nabla G\left( 
\overline{x}\right) =0$, i.e., $\displaystyle\int_{T_{0}}\nabla F\left( t,%
\overline{x}\right) dt^{1}\wedge ...\wedge dt^{p}=0$.

Let us show that 2) implies 3):

The properties of $G$ show that if $\overline{x}$ exists so that $\nabla
G\left( \overline{x}\right) =0$, then $G\left( x\right) \rightarrow \infty $
when $\left| x\right| \rightarrow \infty $, so $\displaystyle%
\int_{T_{0}}F\left( t,x\right) dt^{1}\wedge ...\wedge dt^{p}\rightarrow
\infty $ when $\left| x\right| \rightarrow \infty $.

Now, 3) implies 1). Indeed, the required implication is realized by the
Theorem 1.

{\bf Remark}. The previous results can be extended to PDEs produced in
[12]-[19].

{\bf Aknowledgements.} The authors are grateful to Prof. Dr. Kostake
Teleman and Prof. Dr. Ionel Tevy for their valuable comments on this paper.

\bigskip

\centerline{\bf References}

\bigskip

[1] R. A. Adams: {\it Sobolev Spaces, }Academic Press, 1975.

[2] C. J. Budd, A. R. Humphries, A. J. Wathen: {\it The Finite Element
Approximation of Semilinear Elliptic PDEs in the Cube}, 1997.

[3] T. de Donder: {\it \ Theorie invariantive du calcul des variationes},
1935.

[4] I. Duca, M. Duca: {\it Contractions and Contractive Mappings in
Differential Equations Theory}, U.P.B. Sci. Bull., Series A, Vol. 64, No. 3
(2002), 23-30.

[5 L. V. Kantorovici, G. P. Akilov: {\it Analiz\u{a} func\c{t}ional\u{a}},
Editura \c{S}tiin\c{t}ific\u{a} \c{s}i Enciclopedic\u{a}, Bucure\c{s}ti,
1980.

[6] J. Mawhin, M. Willem: {\it Critical Point Theory and Hamiltonian Systems}%
, Springer-Verlag, 1989.

[7] S. G. Mihlin: {\it Ecua\c{t}ii liniare cu derivate par\c{t}iale},
Editura \c{s}tiin\c{t}ific\u{a} \c{s}i enciclopedic\u{a}, Bucure\c{s}ti,
1983.

[8] R. S. Palais: {\it Foundation of a Global Nonlinear Analysis, }New York,
Benjamin, 1968.

[9] R. E. Showalter: {\it Hilbert Space Methods for Partial Differential
Equations}, Electric Journal of Differential Equations Monograph 01, 1994.

[10] K. Teleman: {\it Metode \c{s}i rezultate \^{i}n geometria diferen\c{t}%
ial\u{a} modern\u{a}}, Editura \c{S}tiin\c{t}ific\u{a} \c{s}i Enciclopedic\u{%
a}, Bucure\c{s}ti, 1979.

[11] C. Udri\c{s}te: {\it From Integral Manifolds and Metrics to Potential
Maps}, Conference Michigan State University, April 13-27, 2001; Eleventh
Midwest Geometry Conference, Wichita State University, April 27-29, 2001,
Atti del' Academia Peloritana dei Pericolanti, Clase 1 di Scienze Fis. Mat.
e Nat., 81-82, A 01006 (2003-2004), 1-14.

[12] C. Udri\c{s}te: {\it Nonclassical Lagrangian Dynamics and Potential Maps%
}, The Conference in Mathematics in Honour of Professor Radu Ro\c{s}ca on
the Occasion of his Ninetieth Birthday, Katholieke University Brussel,
Katholieke University Leuven, Belgium, Dec.11-16, 1999;

http://xxx.lanl.gov/math.DS/0007060.

[13] C. Udri\c{s}te: {\it \ Solutions of DEs and PDEs as Potential Maps
Using First Order Lagrangians}, Centenial Vr\^{a}nceanu, Romanian Academy,
University of Bucharest, June 30-July 4, (2000);

http://xxx.lanl.gov/math.DS/0007061; Balkan Journal of Geometry and Its
Applications 6, 1, 93-108, 2001.

[14] C. Udri\c{s}te, I. Duca: {\it Periodical Solutions of Multi-Time
Hamilton Equations,} Analele Universita\c{t}ii Bucure\c{s}ti, 55, 1 (2005),
179-188.

[15] C. Udri\c{s}te, I. Duca: \ {\it Poisson-gradient \ Dynamical \ Systems
with Bounded Non-linearity}, manuscript.

[16] C. Udri\c{s}te, M. Ferrara, D. Opri\c{s}: {\it Economic Geometric
Dynamics, }Geometry Balkan Press, Bucharest, 2004.

[17] C. Udri\c{s}te, M. Neagu:{\it \ From PDE Systems and Metrics to
Generalized Field Theories}, http://xxx.lanl.gov/abs/math.DG/0101207.

[18] C. Udri\c{s}te, M. Postolache: {\it Atlas of Magnetic Geometric Dynamics%
}, Geometry Balkan Press, Bucharest, 2001.

[19] C. Udri\c{s}te, A-M. Teleman: {\it Hamilton Approaches of Fields Theory}%
, IJMMS, 57 (2004), 3045-3056; ICM Satelite Conference in Algebra and
Related Topics, University of Hong-Kong, 13-18.08.02.

\bigskip

University Politehnica of Bucharest

Department of Mathematics

Splaiul Independen\c tei 313

060042 Bucharest, Romania

udriste@mathem.pub.ro

\end{document}